\documentclass{article}

\usepackage{arxiv}
\usepackage{amsmath}

\usepackage[utf8]{inputenc} % allow utf-8 input
\usepackage[T1]{fontenc}    % use 8-bit T1 fonts
\usepackage{hyperref}       % hyperlinks
\usepackage{url}            % simple URL typesetting
\usepackage{booktabs}       % professional-quality tables
\usepackage{amsfonts}       % blackboard math symbols
\usepackage{nicefrac}       % compact symbols for 1/2, etc.
\usepackage{microtype}      % microtypography
\usepackage{lipsum}
\usepackage{graphicx}
\graphicspath{ {./images/} }

\newcommand{\BEQ}{\begin{equation}}
\newcommand{\EEQ}{\end{equation}}

\title{Musings on the Riemann Hypothesis}

\author{
Ali Nadim \\
Institute of Mathematical Sciences\\
Claremont Graduate University\\
Claremont, CA 91711 USA\\
\texttt{ali.nadim@cgu.edu} \\
}

\begin{document}
\maketitle

\begin{abstract}
We present a few ideas on the Riemann Hypothesis based on properties of analytic functions in the complex plane. In particular, we focus on the real and imaginary parts of the Riemann xi ($\xi$) function whose zeros coincide with those of the zeta ($\zeta$) function within the critical strip. We discuss the forms of the zero contour lines of the two conjugate harmonic functions (the real and imaginary parts of xi) and consider where their intersections could conceivably occur. Those intersections would be the roots of both $\xi$ and $\zeta$ functions. The question of whether a zero could occur away from the critical line becomes equivalent to whether the solutions of a pair of Laplace's equations with well-defined boundary conditions in some semi-infinite strip can possess zero contour lines that intersect within that strip.
\end{abstract}

% keywords can be removed
\keywords{Riemann Zeta Function \and Riemann Hypothesis \and Harmonic Functions}

\section{Introduction}
\label{sec:intro}

Our first encounter with the Riemann zeta function \cite{conrey2003}, $\zeta(s)$, is usually in the form of the infinite series
\begin{equation}\label{eq:series} \zeta(s) = \sum_{n=1}^{\infty} \frac{1}{n^s}\,,\end{equation}
which converges in the complex $s$-plane for $\Re(s)>1$. By analytic continuation, this function can be extended to the remainder of the complex plane and there are common integral representations of $\zeta(s)$ that provide the function in various regions of interest. For example, a representation that holds for $\Re(s)>0$ is given by \cite[cf.~(25.5.3)]{NIST:DLMF}
\begin{equation}\label{eq:intrep} 
\zeta(s)=\frac{1}{(1-2^{1-s}) \Gamma(s)} \int_0^\infty \frac{x^s}{e^x+1}\frac{dx}{x} \,,
\end{equation}
in which $\Gamma(s)$ is the well-known gamma function that generalizes the factorial operation to non-integer values. In Appendix~\ref{appA} we show how this integral representation can be used to evaluate $\zeta(s)$ within the critical strip in terms of two infinite series involving elementary functions.

It is known that in the complex $s$-plane, the zeta function has just one simple pole at $s=1$ and infinitely many zeros. These consist of the trivial zeros that are located along the negative real axis at even negative integer values $(s=-2,-4,-6,\cdots)$, and the nontrivial zeros that are known to be located within the so-called critical strip in the complex plane: $0<\Re(s)<1$. To date, more than ten trillion of these nontrivial roots have been found numerically, all of which have a real part equal to 1/2. Indeed, Riemann speculated that all nontrivial zeros of the zeta function are along the critical line $\Re(s)=1/2$ in the complex plane, but despite overwhelming computational support, no definitive proof of this hypothesis has been found. This is one of the unsolved Millennium Prize Problems \cite{sarnak2004} by the Clay Mathematics Institute; see \url{https://www.claymath.org/millennium-problems/}.

In this brief paper, we discuss a potential approach to proving the Riemann Hypothesis based upon properties of conjugate harmonic functions, specifically the real and imaginary parts of the Riemann xi function, $\xi(s)$, which is related to the zeta function and is known to share the exact same roots as the latter within the critical strip.

\section{The xi Function}

Riemann introduced the so-called xi function, $\xi(s)$, which is related to the zeta function by \cite{Weisstein}
\begin{equation}
\label{eq:xidef} 
\xi(s)=\frac{s(s-1)\Gamma(s/2)}{2\, \pi^{s/2}} \zeta(s) \,. 
\end{equation}
It satisfies the functional relation $\xi(s)=\xi(1-s)$, which implies that $\xi(s)$ is symmetric upon reflection through point $1/2+0\,i$ in the complex $s$-plane. For convenience, we define a new complex variable $z$ which is related to $s$ by 
\[ s=1/2+z \,.\]
When $\xi$ is regarded as a function of this new variable $z$, the functional relation becomes $\xi(z)=\xi(-z)$. That is, $\xi(z)$ is symmetric with respect to reflections through the origin in the $z$-plane. With real and imaginary parts of $z$ denoted by the standard variables $x$ and $y$, i.e., with $z=x+iy$, the Riemann Hypothesis would be that all roots of $\xi(z)$ occur along the $y$-axis (where $x=0$), which is the same as the so-called critical line, $s=1/2$. In terms of the new variables, the critical strip is defined by $-1/2<x<1/2$. The functional relation tells us that if there is a root away from the critical line, e.g., in the first quadrant relative to the origin $z=0$, then there will also be one in the third quadrant. Similarly, if there is a root in the second quadrant, there will also be one in the fourth. Upon examining the odd- or even-ness of the real and imaginary parts of $\xi(z)$ with respect to $x$ and $y$, we will actually find that a root in any quadrant will imply the existence of three other roots in each of the other quadrants.

The function $\xi(z)$ is an entire function in the complex plane and has the Taylor series about $z=0$ (i.e., $s=1/2$) given by \cite{Weisstein}
\[ \xi(z)= \sum_{n=0}^{\infty} a_{2n} z^{2n} \,.\]
Only even powers of $z$ occur in this series. Coefficients $a_{2n}$ are real and can be obtained by computing a complicated integral involving the function
\[ \psi(x)=\sum_{m=1}^{\infty} e^{-m^2\pi x} \,.\]
Explicitly \cite{Weisstein}
\[ a_{2n} = 4 \int_1^\infty \frac{d(x^{3/2}\psi'(x))}{dx} \frac{(\ln\sqrt x)^{2n}}{(2n)!}\frac{dx}{x^{1/4}} \,.\]
The exact values of the coefficients are not needed for our subsequent analysis, but the facts that only even powers of $z$ appear in the Taylor series and the coefficients are real are relevant. 

Consider the first few even powers of $z$:
\begin{align*}
z^2 &= (x^2-y^2) + i (2 xy) \\
z^4 &= (x^4 -6 x^2 y^2+y^4) + i (x^3 y -x y^3) \\
z^6 &= (x^6 -15 x^4 y^2 +15 x^2 y^4 -y^6) + i (6 x^5 y -20 x^3 y^3 +6 x y^5)\,.
\end{align*}
By induction, it is easy to establish that the real parts of all of these and higher even powers of $z$ are even functions in both $x$ and $y$, while the imaginary parts are all odd functions in both $x$ and $y$. Since the coefficients in the series for $\xi(z)$ are themselves real, we find that upon writing $\xi(z)$ in terms of its real ($f$) and imaginary ($g$) parts in the form
\[ \xi(z)=f(x,y)+i g(x,y)\,,\]
function $f(x,y)$ is even in each of its arguments, and $g(x,y)$ is odd in each of its arguments. It becomes clear that the imaginary part $g$ is identically zero along both the $x$- and $y$-axes. Also, the derivative $\partial f / \partial x$ is zero along the $y$-axis, i.e., along the critical line. Also, if there is a root of $\xi$ at some point $(x_o,y_o)$, there will be a root at all four points $(\pm x_o, \pm y_o)$, i.e., in all four quadrants. 

One could draw similar conclusions by replacing $z$ with its polar form, $z=r e^{i\theta}$, in the Taylor series. Each term in even powers of $z$ then involves the factor $e^{i(2n\theta)}$ whose imaginary part is $\sin(2 n \theta)$. Again, along the primary axes (where $\theta=0,\pi/2,\pi,3\pi/2$), the imaginary part of the Taylor series is seen to be zero, with function $\xi$ being purely real along those axes. 

In Appendix~\ref{AppB}, we exploit the fact that multiplying $\zeta(s)$ by the factor shown in Eq.~(\ref{eq:xidef}) produces a function whose imaginary part is identically zero along the critical line, to obtain asymptotic approximations to the argument of the zeta function and the amplitude of the xi function along that line.

\section{Zero Contours of Real and Imaginary Parts}

\begin{figure} 
    \centering
    \includegraphics[scale=1]{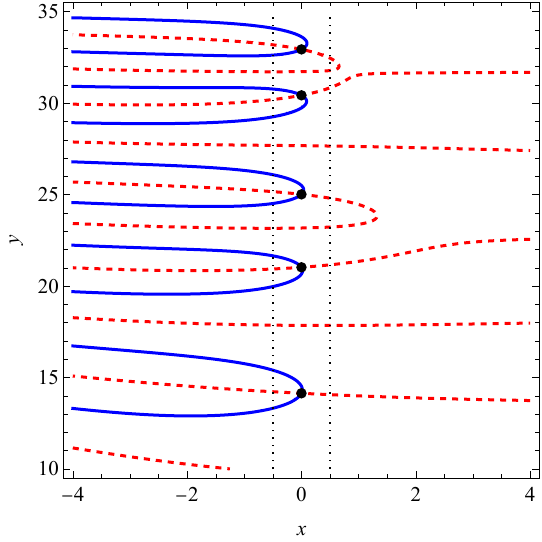}
    \caption{Zero contours of the real part (solid, blue curves) and the imaginary part (dashed, red curves) of $\zeta(z)$. The dotted vertical lines delineate the critical strip. Black dots identify the zeros of the $\zeta$ function in this region. They lie on the critical line $x=0$.}
    \label{fig:ZetaContours} 
\end{figure}

Let us start with the zeta function itself. For $\zeta(z)$ to be zero, both its real and imaginary parts must be zero. In Figure~\ref{fig:ZetaContours}, we have used Mathematica to plot the zero contours of both these functions over some arbitrary domain: $-4<x<4$ and $10<y<35$. The solid blue curves in the figure are the zero contours of the real part of $\zeta$ while the red dashed curves are those of the imaginary part. As seen in the figure, there are multiple intersection points between these lines. These intersection points happen to fall on the critical line (the $y$-axis in this coordinate system). However, it would be impossible to visually determine this. The roots could be slightly away from the axis, and a visual inspection would not reveal that. As seen in the following, this will not be an issue when these contour lines are plotted for the xi function. Where the points of intersection occur, provided that the zeros are simple first-order roots, the zero contours of the real and imaginary parts are perpendicular to each other. This may not seem obvious in the figure, especially since the horizontal and vertical scales are different, but upon zooming in to each root, it can be verified. Theoretically, this would have to be true, since in the vicinity of a first order root $z_o$, the function behaves like $\alpha (z-z_o)$ where $\alpha$ could be a complex coefficient. The zero real and imaginary parts of $z-z_o$ coincide with the horizontal and vertical axes (being perpendicular), and multiplication by $\alpha$ simply rotates these contours by the same angle. 

\begin{figure}
    \centering
    \includegraphics[scale=0.7]{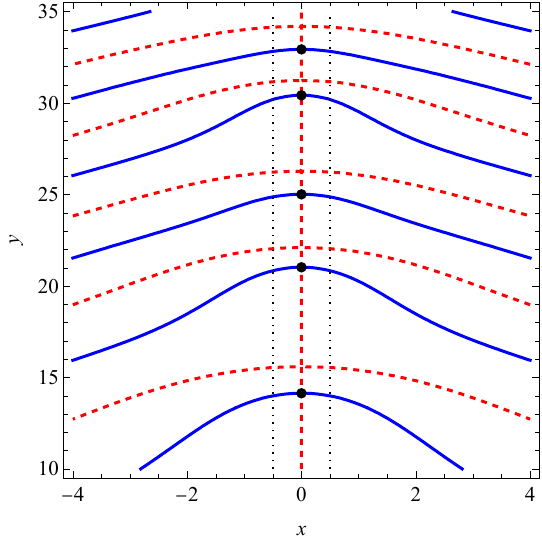}\qquad
    \includegraphics[scale=0.7]{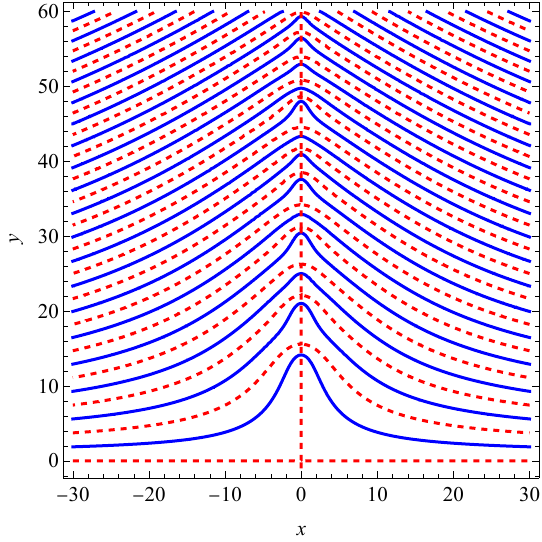}
    \caption{Zero contours of the real part (solid, blue curves) and the imaginary part (dashed, red curves) of $\xi(z)$. The dotted vertical lines delineate the critical strip. Black dots identify the zeros of the $\zeta$ (and $\xi$) function in this region. They lie on the critical line $x=0$. The right plot shows the same contour lines over an extended range.}
    \label{fig:XiContours} 
\end{figure}

We now move on to the plot of the zero contours of the real and imaginary parts of $\xi(z)=f(x,y)+i g(x,y)$, displayed in Figure~\ref{fig:XiContours} over the same range as before (left plot) and an extended range (right plot). In this case, the zero contours $g(x,y)=0$ drawn in red dashed style include the vertical axis, consistent with the fact that the imaginary part of $\xi$ is zero along the critical line. Locally, these blue curves (on which $f=0$) cross the critical line horizontally provided that the corresponding root is a simple first-order zero. Thus, whenever one of the blue curves crosses the vertical axis, we see one of the zeros of both $\xi$ and $\zeta$ along the critical line. There are infinitely many such zeros along the critical line. The question is whether there can be any zeros away from the critical line, i.e., whether intersections of blue and red curves could occur away from the vertical axis in the critical strip. Riemann's Hypothesis would suggest that these cannot occur.

A useful property of these zero contour lines follows from the Cauchy-Riemann (CR) relations governing the real and imaginary parts of analytic functions \cite{carrier2005}. Since $f=0$ along the blue curves, the tangential derivative of $f$ along those curves is also zero: $\partial f/\partial \tau=0$ (we use $\tau$ to denote the tangential direction along these contours). By the CR relations, the normal derivative of $g$ must then also be zero along the same blue curves: $\partial g/\partial n =0$ ($n$ denotes the normal direction at any point along the contour). This means that the blue curves identify a ``ridge line'' for the $g$ surface. Crossing that ridge line in the lateral direction, one passes through a maximum or minimum of $g$. However, since $g$ is a harmonic function and cannot have an actual local maximum or minimum at a point \cite{axler2013}, that ridge line cannot culminate at an actual peak or trough. If we reach a point along the ridge (blue curve) at which the tangential derivative of $g$ is also zero (making the full gradient of $g$ zero there), that point must be a saddle point of the $g$ surface; it cannot be a local minimum or maximum. As such, the zeros along the critical line are saddle points of $g$. The same is true of the dashed red curves, which are zero contours of $g$, which means the tangential derivative of $g$ along them vanishes. As such, the normal derivative of $f$, $\partial f/\partial n$, also vanishes along the red dashed curves, making those ridge lines of the $f$ surface. In particular, along the critical line $x=0$, we must have $\partial f /\partial x=0$, which is consistent with $f$ being an even function of variable $x$. 

Let us denote the value of $f(x,y)$ along the critical line by $F(y)$, i.e., $F(y)=f(0,y)$. Function $F(y)$ has infinitely many zeros along the vertical axis. If those are simple roots, the function must change sign when crossing those points. So between any pair of neighboring roots, $F(y)$ is everywhere either positive or negative, alternating from one segment to the next. Consider traversing up along the $y$ axis from one root to the next, in a segment where $F$ is positive. $F$ must increase from zero, reach at least one maximum along the way, and decrease back down to zero at the next root. At that maximum point of $F(y)$, there must be a saddle point of $f(x,y)$ (since both tangential and normal derivatives of $f$ are zero there) as well as a saddle point of $g(x,y)$ (by the CR relations). These saddle points of $g$ are characterized by the intersection of two red dashed curves along the $y$ axis. In Figure~\ref{fig:XiContours} we see that in the region plotted, there is exactly one such red curve crossing the critical line horizontally in between each neighboring pair of blue curves. In theory, there can be multiple such red curves between neighboring zeros, but that would only occur if the function $F(y)$ passed through multiple points where $F'(y)=0$ in between the two roots; for example $F(y)$ could pass through a maximum, a minimum, and another maximum before returning to zero at the next root, all the while remaining positive on the segment. In that case, we would see three dashed red curves (three saddle points) between the two horizontal blue curves. It is the profile of $F(y)$ along each segment between two neighboring roots that fully determines this structure. 

Let us now consider the question of whether the blue and red curves that cross the critical line horizontally can ever intersect as we move to the right or left, away from the $y$-axis. The answer turns out to be no. Keep in mind that these contours themselves must continue all the way to infinity since, as contour lines of non-constant harmonic functions, they cannot simply terminate somewhere or create closed loops \cite{carrier2005, axler2013}. If any neighboring pair of blue and red curves emanating from the vertical axis ever intersect as we move left or right, the contours would enclose a closed region of the complex plane along whose boundaries either $f=0$ or $g=0$ (i.e., $\partial f/\partial n=0$). Since $f$ is a harmonic function satisfying $\nabla^2 f=0$, its unique solution in the closed region, with boundary conditions that either $f=0$ or $\partial f / \partial n=0$ along the enclosing boundaries, is $f(x,y)=0$ throughout. Its harmonic conjugate function $g(x,y)=0$ will also be zero throughout that closed region, which makes $\xi(z)=0$ in that region. Since $\xi(z)$ is an entire function, by analytic continuation away from this closed region, $\xi(z)$ would have to be identically zero over the entire complex plane, which would be a contradiction. Therefore, the blue and red curves that emanate horizontally from the vertical axis can never intersect with each other or with themselves. This means that no new zeros of the $\zeta$ function can be created by the intersections of the contour lines $f=0$ and $g=0$ that also intersect the critical line. But this leaves open the possibility that additional zero contours of $f$ and/or $g$ that we do not see in Figure~\ref{fig:XiContours} and that reside entirely away from the critical line may occur and may intersect with the current zero contours or with themselves to give rise to off-axis roots of $\xi$. We explore this possibility further in the next section. 

\section{Conformal Mapping onto a Semi-Infinite Rectangular Strip}

\begin{figure}
    \centering
    \includegraphics[trim={3cm 4cm 3cm 4cm}]{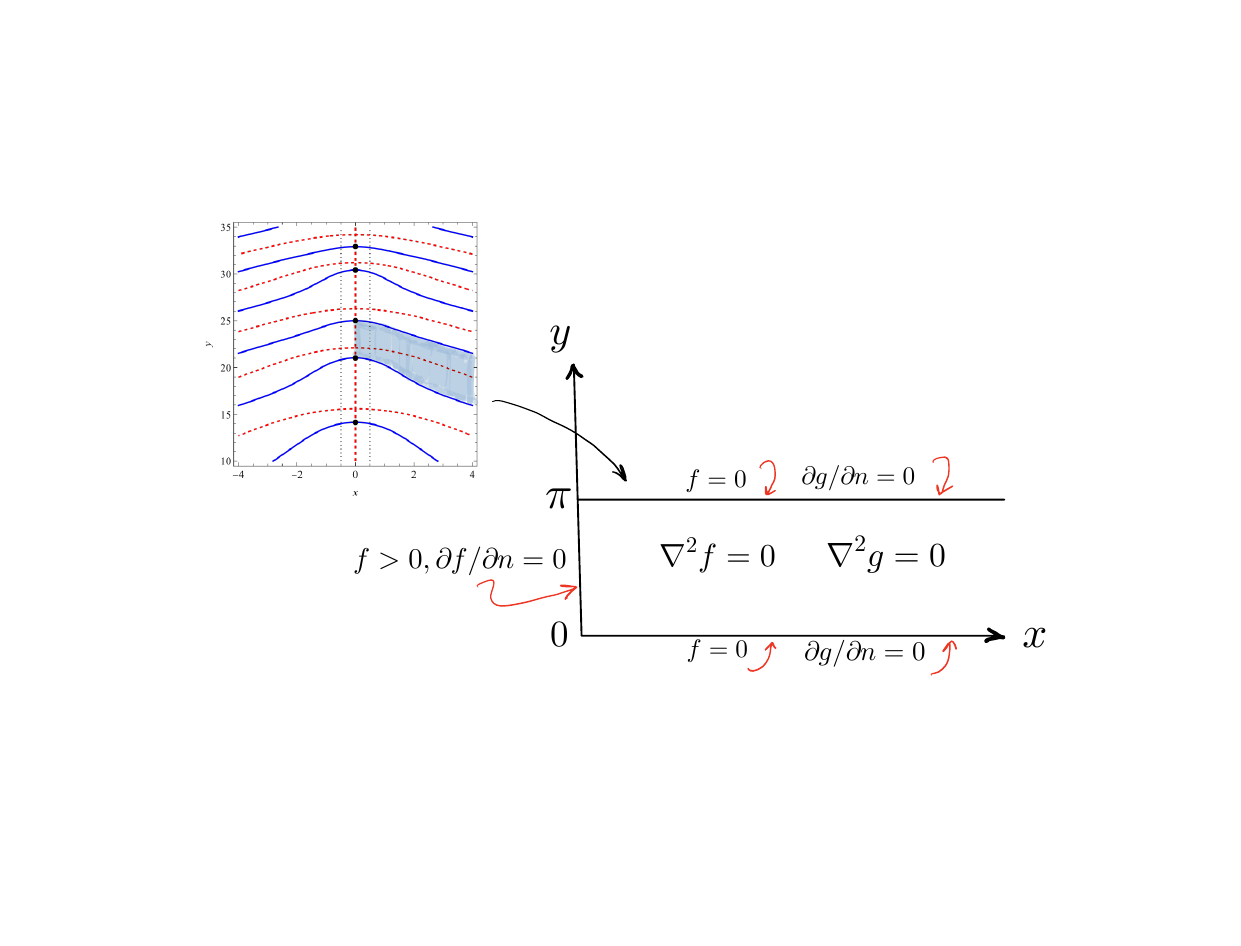}
    \caption{A hypothetical conformal mapping of the shaded region onto the semi-infinite rectangular strip, $0<y<\pi$ and $0<x<\infty$, in a new complex plane.}
    \label{fig:ConformalMap}
\end{figure}

The function $F(y)$ between two neighboring zeros on the critical line determines the structure of zero contours of $f$ and $g$ away from the line. Imagine mapping one of the semi-infinite strips bounded by the critical line on the left (on which $g=0$), and by two successive blue curves emanating from the line on the right side (along which $f=0$), onto a semi-infinite rectangle in a new complex plane, as depicted in Figure~\ref{fig:ConformalMap}. In theory, based on the Riemann Mapping Theorem \cite{carrier2005}, this should always be possible, though we are not able to provide an explicit form for this mapping. For convenience, we still use $x$ and $y$ as the real and imaginary coordinates in the new complex plane, and we take the rectangular region to span $0<y<\pi$ and $0<x<\infty$. Within this rectangular region, it becomes easy to solve for the real and imaginary parts of the mapped function $\xi(x,y)$ given that they are harmonic functions (satisfying Laplace's equation) with well-defined boundary conditions. 

We should clarify that under a conformal map, the function $F(y)$ on the left edge of the original domain gets modified to a new function on the mapped domain. To see this, suppose that we temporarily denote the new coordinates after the mapping by $x_1$ and $y_1$, corresponding to the complex $z_1$-plane. If the conformal mapping is given by analytic function $w$ with real and imaginary parts $u$ and $v$, so that $z_1=w(z)=u(x,y)+iv(x,y)$, the left edge of the original domain $x=0$ gets mapped to the left edge $x_1=u(0,y)=0$, and the corresponding $y_1$ coordinate will be $y_1=v(0,y)$ along the new left edge, chosen to be in the range $0<y_1<\pi$ as $y$ varies between two neighboring zeros. The function $F(y)$ therefore gets mapped to $\tilde{F}(y_1)$, with $y_1$ and $y$ being related by $y_1=v(0,y)$. The inverse of this function would be needed to obtain $y(y_1)$ along the edge, using which $\tilde{F}(y_1)=F(y(y_1))$. However, for simplicity of notation, we still use $(x,y)$ in the mapped domain and refer to the function on the left edge by $F(y)$.

To obtain the real part $f(x,y)$, we need to solve the following boundary value problem:
\[ \nabla^2 f = 0\,\quad (0<x<\infty\,, \quad 0<y<\pi)\]
subject to the boundary conditions:
\begin{align*}
f = 0 & \quad\mbox{at}\quad y=0\,,~ y=\pi \\
f = F(y)\,,~ \partial f/\partial x = 0 & \quad\mbox{at}\quad x=0 \,.
\end{align*}
No boundary condition is imposed at infinity, i.e., as $x\rightarrow\infty$ for $0<y<\pi$. This makes the problem a little different from those typically encountered. To solve such an elliptic partial differential equation, it is more customary to have the boundary condition that $f$ goes to zero at infinity, rather than impose two boundary conditions on the same left edge of the domain. However, since $f(x,y)$ needs to be an even function of $x$, its $x$-derivative must vanish at $x=0$, and we know the value of $f$ on the left boundary as well, given by $F(y)$. So imposing both boundary conditions is inevitable. Finding the solution to this problem is still straightforward using separation of variables although $f(x,y)$ grows (rather than decays) exponentially at infinity. Indeed, the function $\xi(z)$, although entire, must have an essential singularity at infinity \cite{carrier2005}, which is consistent with this solution. 

By separation of variables, the general solution to $f(x,y)$ is given by
\[ f(x,y)=\sum_{n=1}^\infty b_n \sin(ny) \cosh(nx) \,.\]
The coefficients $b_n$ can be found by satisfying the boundary condition at the left edge: $F(y)=\sum_{n=1}^\infty b_n \sin(ny)$, from which, due to the orthogonality of the sine functions, one obtains
\[ b_n = \frac{2}{\pi} \int_0^\pi F(y) \sin(ny) \,dy \,.\]
Assuming $F(y)>0$ on the segment $0<y<\pi$, the leading coefficient $b_1$ must be positive. Since we are primarily interested in finding the zeros of $f(x,y)$ within the strip, we can always factor out this coefficient $b_1$, with the modified series having a coefficient of 1 in front of the first term, and the rest of the coefficients being replaced by their ratios to $b_1$. This is equivalent to setting $b_1=1$ without loss of generality.

The conjugate function $g(x,y)$ can be obtained from $f(x,y)$ by integrating the Cauchy-Riemann equations: $\partial f/\partial x = \partial g / \partial y$, and $\partial f/\partial y = -\partial g / \partial x$. Alternatively, it can be obtained independently by using separation of variables to solve
\[ \nabla^2 g = 0\,\quad (0<x<\infty\,, \quad 0<y<\pi)\]
within the strip, with boundary conditions:
\begin{align*}
\partial g /\partial y = 0 & \quad\mbox{at}\quad y=0\,,~ y=\pi \\
g=0\,,~ \partial g/\partial x = -F'(y) & \quad\mbox{at}\quad x=0 \,.
\end{align*}
Note that the normal derivative of $g$ along the left boundary is being related to the tangential derivative of $f$ on that boundary via the CR relations. The solution is found to be
\[ g(x,y)=-\sum_{n=1}^\infty b_n \cos(ny) \sinh(nx) \,.\]
As before, we can factor out $b_1$ or set it to unity and work with the modified series.

If $F(y)$ [actually $\tilde{F}(y_1)$] is exactly proportional to $\sin(y)$ leaving only the leading term in the series solution, the exact solutions to $f$ and $g$ turn out to be
\[ f(x,y)=\sin(y)\cosh(x) \,,\quad g(x,y)=-\cos(y)\sinh(x)\,.\]
In that case, we see that $f$ does not have any new zeros (or, more precisely, zero contours) within the semi-infinite strip. It starts positive along the $y$ axis and increases from there as $x$ goes to infinity. The conjugate function $g(x,y)$ has a zero contour along the horizontal line $y=\pi/2$. No new zeros of $\xi(z)$ occur within the strip. Therefore, with the boundary condition $F(y)=\sin(y)$, it would be impossible to encounter an off-axis zero. If we were lucky enough that the conformal mapping always produced this $\tilde{F}(y_1)$ on the left edge of the mapped domain, the Riemann Hypothesis would be proven. However, $F(y)$ is not expected to be that simple, so more general solutions may include additional, if not all, terms of the series obtained by separation of variables. 

\begin{figure}
    \centering
    \includegraphics[scale=0.5]{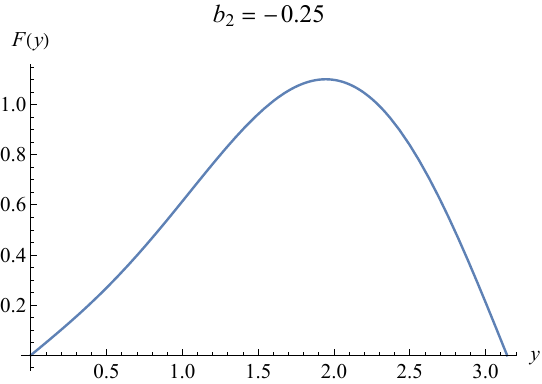} \hspace{0.9cm} 
    \includegraphics[scale=0.5]{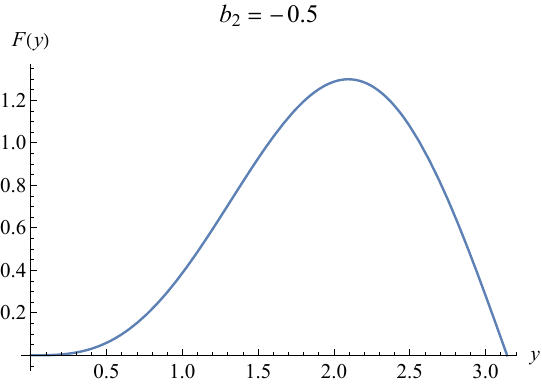} \hspace{0.9cm}
    \includegraphics[scale=0.5]{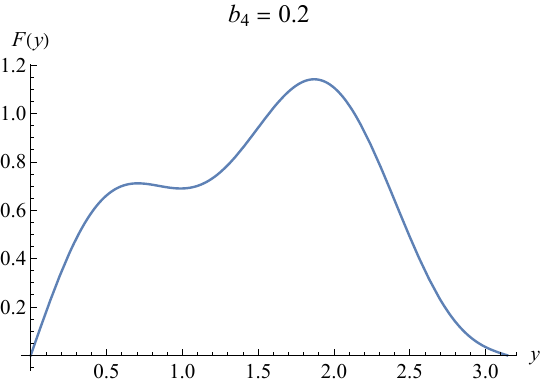}
    \\[0.5cm]
    \includegraphics[scale=0.55]{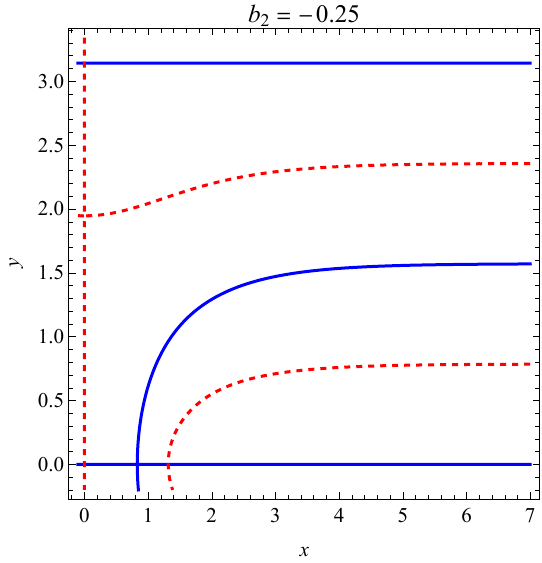} \hspace{0.5cm} 
    \includegraphics[scale=0.55]{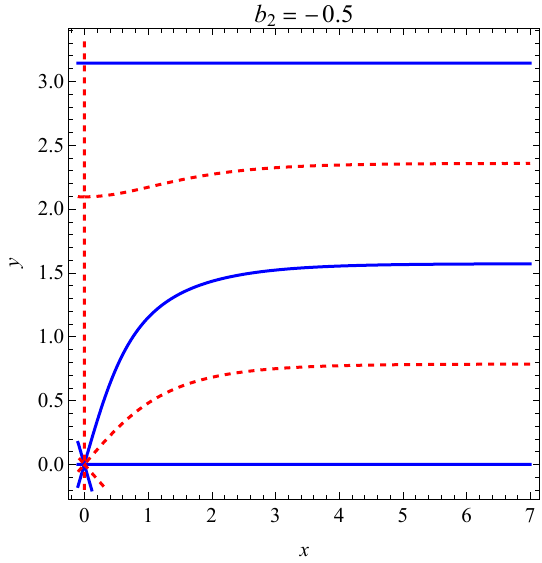} \hspace{0.5cm}
    \includegraphics[scale=0.55]{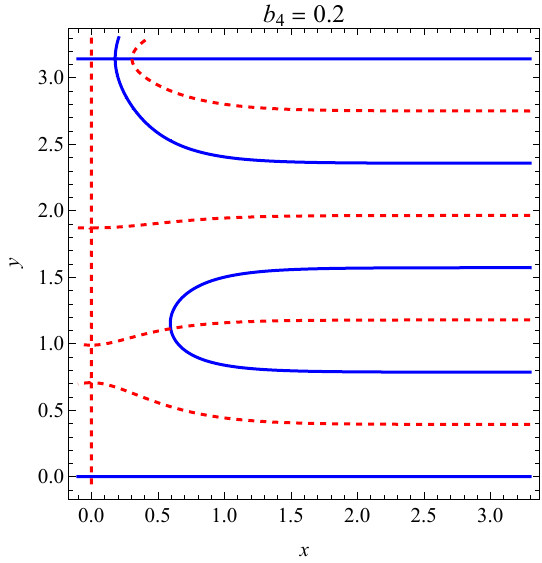}
\caption{Solution via separation of variables. The top graphs show the profile $F(y)$ along the left edge of the semi-infinite strip, each with just one extra mode in its Fourier sine series; the bottom plots show the corresponding zero contours of $f$ (solid blue curves) and $g$ (dashed red curves). Any intersection between the red and blue curves is a zero of the $\xi=f+i\,g$.}
\label{fig:SepVarSoln}
\end{figure}

To explore the question of whether off-axis zeros are even possible, we choose to examine a few simple cases where only one of the next several terms is included in the solution. For instance, corresponding to 
\[ F(y)=\sin(y)+b_2 \sin(2y)\,,\]
we would have the solutions
\begin{align*}
    f(x,y) &= \sin(y) \cosh(x) + b_2 \sin(2y) \cosh(2x) \,,\\
    g(x,y) &= -\cos(y) \sinh(x) - b_2 \cos(2y) \sinh(2x) \,.
\end{align*}
The additional mode $\sin(2y)$ in $F(y)$ is positive over $0<y<\pi/2$, but negative over $\pi/2<y<\pi$. Moreover the exponentially growing term $\cosh(2x)$ grows faster than $\cosh(x)$ as $x$ increases. So no matter how small the coefficient $b_2\neq0$ might be, the second term in the solution to $f(x,y)$ dominates for large enough $x$. So in this case, $f(x,y)$ will develop a new zero contour line within the semi-infinite strip as $x$ gets large. This is also true if any of the other higher modes were included in the solution. Now both $f$ and $g$ will have zero contours within the strip and the question becomes whether these contours can intersect and create a new zero for the $\xi$ function away from the critical line. When we plot the zero contours of $f$ and $g$ based on this solution over the ranges of acceptable $b_2$ values---to keep $F(y)>0$, we must restrict $b_2$ to the range $-0.5<b_2<0.5$---we find that new off-axis zeros can occur. As $b_2$ approaches either 0.5 or $-0.5$, one of original roots along the critical line becomes a double root. The left and middle panels of Figure~\ref{fig:SepVarSoln} show the solutions for $b_2=-0.25$ and $b_2=-0.5$. The top graphs show $F(y)$ versus $y$ and the bottom graphs show the corresponding zero contours of $f$ (solid blue curves) and $g$ (dashed red curves). We see that for $b_2=-0.25$ there are new zero contours of both $f$ and $g$. They do not intersect each other, but the new red dashed curve on which $g=0$ does intersect the existing $f=0$ line at the bottom edge of the domain, thus creating a new zero of the complex function $f+ig$. As $b_2$ approaches $-0.5$, these new zero contours of $f$ and $g$ merge with the previous zero along the left edge and make that into a double zero. As can be seen in the middle panels of Figure~\ref{fig:SepVarSoln}, the function $F(y)$ becomes tangent to the horizontal line at $y=0$ and the adjacent contour lines $f=0$ and $g=0$ are no longer perpendicular and instead divide the complex plane into sectors of equal angle $\pi/6$ locally near the bottom left corner at the double root.

As an additional example, we also consider the case where $F(y)=\sin(y)+b_4 \sin(4y)$. When $b_4=0.2$, the function $F(y)$ becomes wavy on the segment and has two maxima and one minimum on $0<y<\pi$ between the two end points where it is zero. These are shown on the right panels of Figure~\ref{fig:SepVarSoln}. In the plot of the zero contour lines of $f$ and $g$ shown on the bottom right, we see that three red dashed curves emanate from the left edge, one at each location where $F'(y)=0$. In addition to those, new zero contours of $f$ and $g$ also form, with intersections both within the interior of the strip and at the top boundary, giving rise to two new off-axis zeros.

These examples illustrate that simply having $F(y)$ be positive or negative over the entire segment between two zeros is not sufficient to avoid new zeros of the complex function $f+ig$ away from the critical axis. Obviously if $F(y)$ [i.e., $\tilde{F}(y_1)$] happens to be exactly proportional to $\pm\sin(y)$, then no zeros can be found away from the critical line, but this is unlikely. However, if $F(y)$ is not simply expressible by a finite truncation of a Fourier sine series, but has infinitely many terms, or if one can find criteria on $F(y)$ that would make it impossible for new zero contours of $f$ and $g$ to form and intersect with each other or with the previous zero contours, then some progress on the Riemann Hypothesis could be made. The $F(y)$ on each segment between neighboring zeros needs to satisfy those criteria, which may not be easy to show. 

In the next section, we use an approach based on Taylor expansion near the critical axis to show that if the successive even-order derivatives of $F(y)$ alternate in sign, then $f$ cannot develop new zero contours away from the critical line. This is certainly true for $F(y)=\sin(y)$, which changes sign each time it is differentiated twice. But it is not clear if the $F(y)$ that arises from evaluating the real part of the xi function along the critical axis, or its counterpart after conformal mapping, would satisfy this criterion.

\section{Taylor Series near the Critical Line}

In searching for zeros of the xi (and zeta) function away from, but close to, the critical line $x=0$, some conclusions can also be drawn by considering the Taylor series for $\xi(z)=f(x,y)+ig(x,y)$ near the vertical axis $x=0$ at some arbitrary $y$ location. For the real part $f(x,y)$, the Taylor series about $x=0$ reads:
\[ f(x,y)=f(0,y)+x \frac{\partial f}{\partial x}(0,y)+\frac{1}{2}x^2 \frac{\partial^2 f}{\partial x^2}(0,y)+\cdots\,.\]
Upon denoting $f(0,y)=F(y)$ as we had done earlier and recognizing that $\partial f/\partial x=\partial g/\partial y=0$ along the line $x=0$ (based on the CR relations), and using the fact that $\nabla^2f=0$ which implies that $\partial^2f/\partial x^2=-F''(y)$ along the axis, the first two nonzero terms in the Taylor series become
\begin{equation}
f(x,y)=F(y)-\frac{1}{2}x^2 F''(y)+\cdots\,.    
\label{eq:parabola}
\end{equation}
The third derivative $\partial^3f/\partial x^3$ along $x=0$ can be shown to be zero by combining the CR relations and Laplace's equation for $f$. At a point where $F(y)$ is positive, if $F''(y)$ is negative, the parabola in $x$ represented by Eq.~(\ref{eq:parabola}) increases as we move away from the critical axis, which means that higher order terms would be needed if the function were to turn around and reach zero even farther away. Where $F(y)$ starts negative, if $F''(y)$ is positive, the same type of conclusion can be drawn. The next nonzero term in the Taylor series will be the $x^4$-term. Indeed, only even powers of $x$ can appear since $f$ must be an even function of $x$. Its explicit calculation involves simplifying the 4th derivative $\partial^4f/\partial x^4$, which can be done by taking the second partial derivative of $\nabla^2f=0$ with respect to $x$ to yield
\[ \frac{\partial^4f}{\partial x^4}=-\frac{\partial^2}{\partial x^2}\frac{\partial^2}{\partial y^2}f=-\frac{\partial^2}{\partial y^2}\frac{\partial^2}{\partial x^2}f=\frac{\partial^4f}{\partial y^4}\,.\]
When evaluated at $x=0$, the latter is the same as ${d^4F}/{dy^4}\equiv F^{(4)}(y)$. Higher order terms can be simplified in a similar manner and the extended Taylor series reads
\begin{equation} f(x,y)=F(y)-\frac{1}{2}x^2 F''(y)+\frac{1}{4!}x^4 F^{(4)}(y)-\frac{1}{6!}x^6 F^{(6)}(y) +\cdots\,.
\label{eq:extendedTaylor}
\end{equation}
Eq.~(\ref{eq:extendedTaylor}) suggests that if the even-order derivatives of $F(y)$ (i.e., $F(y)$, $F''(y)$, $F^{(4)}(y)$, $F^{(6)}(y)$, \ldots) alternate in sign (some could be zero), there can be no new zeros of $f(x,y)$ horizontally away from the $y$-axis at $y$-positions where $F$ starts out positive or negative. If $F(y)$ resulting from evaluating $\xi$ along the critical axis (or $\tilde{F}(y_1)$ after the conformal mapping) had this property everywhere, the Riemann Hypothesis would hold. This only constitutes a sufficient condition, however, and not a necessary one, to preclude off-axis zeros at a particular $y$ value.

Similarly, one can obtain the Taylor series for $g(x,y)$ near the axis $x=0$ along which $g=0$. We have
\[ g(x,y)=g(0,y)+x \frac{\partial g}{\partial x}(0,y)+\frac{1}{2}x^2 \frac{\partial^2 g}{\partial x^2}(0,y)+\frac{1}{6}x^3 \frac{\partial^3 g}{\partial x^3}(0,y)+\cdots\,.\]
Using the CR relations and the fact that $\nabla^2 g =0$, the first two nonzero terms that remain are given by:
\[ g(x,y)=-xF'(y)+\frac{x^3}{6}F'''(y)+\cdots\,.\]
Additional terms in this Taylor series can also be found, involving odd powers of $x$ and higher-order odd derivatives of $F(y)$ with alternating signs. For example, the next nonzero term in the Taylor series would be $-\frac{1}{5!}x^5 F^{(5)}(y)$. A similar conclusion regarding the impossibility of new zero contours of $g(x,y)$ away from the critical axis can be drawn, starting in regions where $F'(y)\neq 0$, when the odd-order derivatives of $F$ alternate in sign.

It should be noted that the function $F(y)$ in the above Taylor series cannot be any arbitrary function. It must come from the evaluation of a harmonic function $f(x,y)$ along the vertical axis $x=0$, one whose normal derivative $\partial f/\partial x$ along that axis is zero, corresponding to its harmonic conjugate $g(x,y)$ being zero there as well. As an example, if we take $F(y)=(\pi/8)y(\pi-y)$ (the factor $\pi/8$ makes $b_1=1$ in its Fourier sine series), a corresponding harmonic function $f(x,y)$ that reduces to $F(y)$ at $x=0$ could be $f(x,y)=F(y)+\pi x^2/8$. This is consistent with the above Taylor series, terminating after the first two nonzero terms. Function $f(x,y)$ does not become zero away from the vertical axis $x=0$, in the range $0<y<\pi$ over which $F(y)>0$. This matches the expectation based on the criteria developed above. (Zero contours of this function do exist, but they are above and below the strip $0<y<\pi$.) However, if we were solving $\nabla^2f=0$ in the semi-infinite rectangular strip with this particular $F(y)$ along the left edge and with the additional boundary conditions that $f=0$ at $y=0$ and $y=\pi$, we would get a different harmonic function using separation of variables, namely:
\[ f(x,y) = \sum_{k=0}^{\infty} \frac{1}{(2k+1)^3} \sin[(2k+1)y]\cosh[(2k+1)x] \,.\]
At $x=0$ this reduces to the same $F(y)=(\pi/8)y(\pi-y)$ function (represented by a Fourier sine series). In this case when its zero contour lines are plotted together with the zero contour lines of its harmonic conjugate function $g$ (after a finite truncation of both series), we do see intersecting zero contour lines in the strip and therefore off-axis zeros (plots not shown). However, since the higher terms in such series dominate for large $x$, it is not clear that a finite truncation of the infinite series is even appropriate.

\begin{figure}
    \centering
    \includegraphics[scale=0.55]{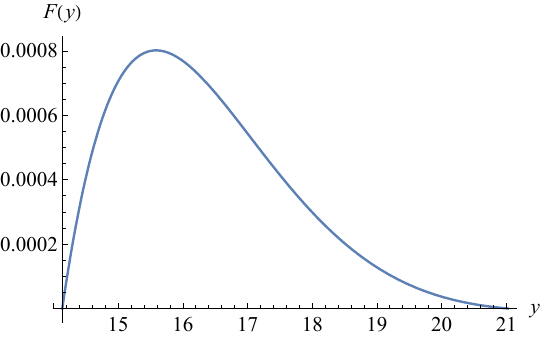} \hspace{0.5cm} 
    \includegraphics[scale=0.55]{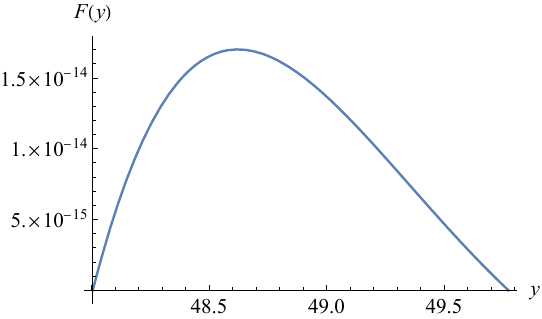} \hspace{0.5cm}
    \includegraphics[scale=0.55]{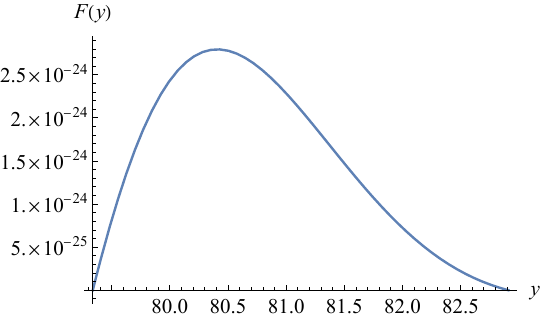}
\caption{Absolute value of the real part of the xi function between neighboring roots along the critical line: $|F(y)|$ vs.~$y$; left plot: $14.1347<y<21.022$, middle plot: $48.0052<y<49.7738$, right plot: $79.3374<y<82.9104$ .}
\label{fig:SampleFofYs}
\end{figure}

For reference, in Figure~\ref{fig:SampleFofYs}, we show the true $F(y)$ function, taken to be the absolute value of the real part of the $\xi$ function along the critical line, between three pairs of neighboring roots (1st and 2nd, 9th and 10th, 21st and 22nd) to give the reader some sense for what these look like. Note the vertical scale in these plots. The overall magnitude of the $\xi$ function decreases exponentially fast as we traverse up the critical axis, though it continues to oscillate and cross zero at successive roots. As evident, these functions themselves do not appear to have negative second derivatives over the entire interval. It is unknown whether the conformal mapping that maps each semi-infinite strip onto the semi-infinite rectangle (as in the last section) might modify these $F(y)$ functions onto mapped functions $\tilde{F}(y_1)$ that are either proportional to $\sin(y_1)$ or satisfy the sufficient conditions (alternating signs of even or odd derivatives) that preclude off-axis zeros. 

\section{Conclusions}

We have presented a few basic ideas in the hope that they might inspire new approaches to proving the Riemann Hypothesis. Mainly, since the shape of the real part of the Riemann xi function, $F(y)=f(0,y)$, along the critical line fully determines the function away from that line [the imaginary part being zero there: $g(0,y)=0$], an idea would be to find the criteria on $F(y)$ that preclude additional zero contour lines of $f(x,y)$ and $g(x,y)$ to form and/or intersect away from the critical line. One can then see if the given $F(y)$ satisfies those criteria, in which case RH would be established. Conveniently, this may be done segment by segment, since in between neighboring roots of $F(y)$ along the critical line, $F(y)$ is either positive or negative, maintaining its sign, and the zero contour lines of $f(x,y)$ that emanate from those roots extend to infinity on either side without intersecting, naturally creating semi-infinite strips within which one could theoretically obtain $f$ and $g$ by solving Laplace's equation for each, subject to well-defined conditions on the boundaries of the strip. Based on a theoretical conformal mapping that maps those semi-infinite strips to a rectangular one, we explored the possible forms of the solutions to $f(x,y)$ and $g(x,y)$. We found that if after the conformal mapping, $F(y)$ happens to become proportional to $\pm \sin(y)$ on the mapped interval $0<y<\pi$, then no zeros of $\xi(z)$ could occur away from the critical line. However, it would require great fortune for the conformal mapping that maps the original semi-infinite strip to a rectangular one to always map the original $F(y)$ shape to a $\sin(y)$ profile. Finding the conformal mapping for transforming the strips is by no means easy or trivial, although it does exist based on the Riemann Mapping Theorem. We also examined the Taylor series for $f$ and $g$ in the vicinity of the critical line, showing that both functions can be expressed in terms of $F(y)$ and its derivatives. The series for $f$ involves only even derivatives of $F(y)$ while that for $g$ only the odd derivatives. Because of the alternating pattern of signs in those series, we could come up with simple criteria that preclude additional zeros of $f$ or $g$ away from the critical line. If $F(y)$ is proportional to $\sin(y)$, those criteria are again satisfied, but for the actual $F(y)$ function, they do not appear to hold. Nevertheless, in searching for potential zeros within the critical strip itself, where $-1/2<x<1/2$, an approximate truncated version of the Taylor series may be helpful.

\appendix

\section{Evaluating the zeta Function in the Critical Strip}
\label{appA}

Starting with the integral representation (\ref{eq:intrep}) we can actually obtain an expression involving elementary functions to evaluate $\zeta(s)$ explicitly within the critical strip. In this appendix, let us use the original $s$ variable and refer to its real and imaginary parts by $\sigma$ and $t$, i.e., $s=\sigma + i \,t$, which is common in the zeta function literature. 

In the region $0<x<\infty$ of the integral in (\ref{eq:intrep}), $e^x>1$, so we can expand
\[ \frac{1}{e^x+1}=\frac{1}{e^x}\left(\frac{1}{1+e^{-x}}\right)=e^{-x}\sum_{m=0}^\infty (-1)^m e^{-mx}=\sum_{m=0}^\infty (-1)^m e^{-(m+1)x}\,,
\]
by means of the geometric series. Integrating each term in this series, as needed in Eq.~(\ref{eq:intrep}), involves 
\[ \int_0^\infty x^{s-1} e^{-(m+1)x} dx = (m+1)^{-s} \int_0^\infty u^{s-1} e^{-u} du = \frac{\Gamma(s)}{(m+1)^s}\,,\]
in which substitution $(m+1)x=u$ was made to simplify the integral.
Evaluation of Eq.~(\ref{eq:intrep}) then yields
\begin{equation}\label{eq:zetainstrip}
\zeta(s)=\frac{1}{(1-2^{1-s})}\sum_{m=1}^\infty \frac{(-1)^{m-1}}{m^s}\,.
\end{equation}
Note that the sum in the above is very similar to the original definition of the zeta function except for the alternating signs. It is known at the Dirichlet eta function or the alternating zeta function. Using $s=\sigma+i\,t$, one can write $m^{-s}=m^{-\sigma}m^{-it}=m^{-\sigma}[\cos(t\ln m)-i\,\sin(t \ln m)]$. The denominator of the factor in front of the sum in Eq.~(\ref{eq:zetainstrip}) is written similarly as
\[ 1-2^{1-s}=1-2^{1-\sigma}\cos(t\ln 2) +i\,2^{1-\sigma} \sin(t \ln 2)\,.\]
As $t$ varies, this equation describes a circle of radius $2^{1-\sigma}$ centered at 1 in the complex $s$-plane. Provided $0<\sigma<1$ which makes $2^{1-\sigma}>1$, this factor is always finite and does not pass through zero, so its reciprocal is also never zero. (However, note that for $s=1$, this term becomes zero, which is an indication of the pole in the zeta function at $s=1$.) Combining these results provides the following explicit formula for the zeta function in the critical strip:
\begin{equation}
\label{eq:explicit}
\zeta(\sigma,t)=\frac{\tilde{f}(\sigma,t)+ i\,\tilde{g}(\sigma,t)}{1-2^{1-\sigma}\cos(t\ln 2) +i\,2^{1-\sigma} \sin(t \ln 2)}\,,
\end{equation}
in which
\begin{align} 
\tilde{f}(\sigma,t)& =\sum_{m=1}^\infty\frac{(-1)^{m-1}}{m^\sigma}\cos[t\ln(m)]\,, \\
\tilde{g}(\sigma,t)& =\sum_{m=1}^\infty\frac{(-1)^{m}}{m^\sigma}\sin[t\ln(m)]\,.
\end{align}
The real and imaginary parts of $\zeta$ can be obtained easily from equation (\ref{eq:explicit}) by multiplying its numerator and denominator by the complex conjugate of its denominator. However, to look for the zeros of the zeta function in the strip, we can simply search for the zeros of $\tilde{f}$ and $\tilde{g}$ since the factor multiplying them does not become zero in the critical strip. 

\begin{figure}
    \centering
    \includegraphics[scale=1]{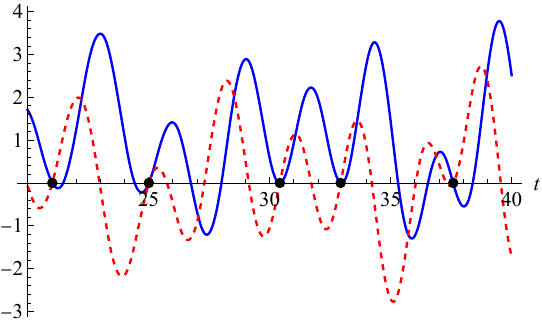}
    \caption{Plots of the real (solid blue) and imaginary (red dashed) parts of the Dirichlet eta function along the critical line $\sigma=1/2$ for $20<t<40$. Where both of these are zero one finds a zero of the zeta function identified by large black dots.}
    \label{fig:etafxn}
\end{figure}

In Figure~\ref{fig:etafxn}, we have plotted the functions $\tilde{f}(\frac{1}{2},t)$ and $\tilde{g}(\frac{1}{2},t)$ versus $t$ by truncating their respective series after 5000 terms and evaluating them over the range of $t$ indicated. We see that both the real and imaginary parts are highly oscillatory functions, although obviously the oscillations are not regular or periodic. Each of the functions has many zeros, and from time to time, both of them are zero at the same $t$ value, which are the nontrivial roots of the zeta function.
An interesting feature of these two functions is that if we think of $t$ as `time' and identify the `frequencies' $\omega$ in the time-dependent components $\sin(\omega t)$ and $\cos(\omega t)$, we see that the set of frequencies in these time series includes the natural logarithms of all integers. Upon writing those integers in terms of their prime number factors, we see that the frequencies involve the logarithms of all the prime numbers, all of their harmonics, and all possible sums of the latter. For example, if $m=p_1^k p_2^\ell$, with prime numbers $p_1$ and $p_2$, and $k$ and $\ell$ their integer powers, frequency $\omega=\ln(m)=k\ln(p_1)+\ell \ln(p_2)$, being a sum of harmonics of the prime frequencies $\ln(p_1)$ and $\ln(p_2)$. As a time series, these functions of $t$ seem to contain all such frequencies, with the result being a very non-periodic function of time containing frequencies that are natural logarithms of all prime numbers and their harmonics and additive combinations.   

\section{Asymptotics along the Critical Line}
\label{AppB}

We observe from Eq.~(\ref{eq:xidef}) that upon multiplying the zeta function $\zeta(s)$ by the factor
\begin{equation}
h(s) = \frac{s(s-1)\Gamma(s/2)}{2\, \pi^{s/2}} \,,
\end{equation}
we produce a new function $\xi(s)$ whose imaginary part is identically zero along the critical line $s=1/2$. The real and imaginary parts of the zeta function itself have isolated zeros along that line but are not identically zero. This implies that if we express the zeta function along the critical line in polar form: $\zeta(t)=r(t)e^{i \theta(t)}$ and do the same for the above factor $h(t)=\rho(t)e^{i \phi(t)}$, with $s=1/2+i\,t$, the two angles $\theta(t)$ and $\phi(t)$ must satisfy the relation: $\theta(t)+\phi(t)=n\pi$ with $n$ being an integer. For the function $h(s)$, it is possible to obtain an asymptotic approximation for large $t$. With the aid of Mathematica, we find that for large $t$,
\begin{equation}
\label{eq:rho}
|h(t)|=\rho(t) \sim \frac{\sqrt{\pi}\, t^{7/4}}{2^{1/4}} e^{-\pi t/4} \quad \mbox{as} \quad t \rightarrow \infty \,,   
\end{equation}
and
\begin{equation}
\phi(t) \sim \frac{7 \pi}{8}-\frac{1}{2}t[ 1+ \ln(2\pi/t)] \sim -t \ln(\sqrt{2\pi e/t})  \quad \mbox{as} \quad t \rightarrow \infty \,.   
\end{equation}
If we plot $|\xi(t)|$ versus $t$ along the critical line (which gives us the $F(y)$ function used in the main part of the paper) but first divide it by $\rho(t)$ from Eq.~(\ref{eq:rho}) (plot not shown), we find that this normalized function remains of order unity rather than decay for large $t$, and on each segment between neighboring zeros, the resulting function (essentially $r(t)=|\zeta(t)|$) looks much more symmetric about the midpoint of the segment and closer to a segment-wise sine function. Also, plotting $\arg[\zeta(t)]+\phi(t)$ versus $t$ for large $t$ (plot not shown), we see a series of horizontal line segments (i.e., piecewise constants) with jump discontinuities of magnitude $\pi$ (sometimes $-\pi$ or $2\pi$), consistent with having $\theta(t)+\phi(t)=n\pi$. It should be noted that the argument function in Mathematica always returns a value between $-\pi$ and $\pi$ while $\phi(t)$ is a continuous function so some of the observed jump discontinuities result from that. 

\bibliographystyle{unsrt}  
\bibliography{references}  

\end{document}